\documentclass{gtart}

\def\ifplaintex{\expandafter\ifx\csname documentclass\endcsname\relax}

\def\gtp{{\mathsurround=0pt\it $\cal G\mskip-2mu$eometry \&\ 
$\cal T\!\!$opology $\cal P\!$ublications}}  

\def\recd{{\small Received:\qua\receiveddate\ifx\reviseddate\relax
\else\qquad Revised:\qua\reviseddate\fi\par}} 

\def\lognumber#1{\def\thelognumber{#1}}
\def\volumenumber#1{\def\thevolumenumber{#1}}
\def\volumeyear#1{\def\thevolumeyear{#1}}
\def\papernumber#1{\def\thepapernumber{#1}}
\def\pagenumbers#1#2{\def\startpage{#1}\def\finishpage{#2}}
\def\published#1{\def\publishdate{#1}}

\def\received#1{\def\receiveddate{#1}}

\def\accepted#1{\def\accepteddate{#1}}
\def\asciititle#1{\def\theasciititle{#1}}

\long\def\asciiabstract#1{\long\def\theasciiabstract{#1}}

\let\\\par\let\thelognumber\relax\let\thevolumenumber\relax
\let\thepapernumber\relax\let\thevolumeyear\relax\let\startpage\relax
\let\finishpage\relax\let\publishdate\relax\let\receiveddate\relax
\let\reviseddate\relax\let\accepteddate\relax\let\theasciititle\relax
\let\theasciiauthors\relax
\let\theasciiabstract\relax

\let\theasciiemail\relax

\ifplaintex
\font\logobig=cmssbx10 scaled 3836
\font\logomed=cmssbx10 scaled 2557
\else
\font\logobig=cmssbx10 scaled 4200
\font\logomed=cmssbx10 scaled 2800
\fi

\long\def\makeagttitle{   
\count0=\startpage
\agt\hfill      
\hbox to 45truept{\vbox to 0pt{\vglue -13truept{\logomed A\kern -.37em{\logobig 
T}\kern -.38em G}\vss}\hss}
\break
{\small Volume \thevolumenumber\ (\thevolumeyear)
\startpage--\finishpage\nl
Published: \publishdate}

\vglue .25truein

{\parskip=0pt\leftskip 0pt plus
1fil\def\\{\par\smallskip}{\Large\bf\thetitle}\par\medskip} \vglue
0.05truein

{\parskip=0pt\leftskip 0pt plus 1fil\def\\{\par}{\sc\theauthors}
\par\medskip}%
 
\vglue 0.03truein 

{\small\leftskip 25truept\rightskip 25truept{\bf Abstract}\stdspace\theabstract

{\bf AMS Classification}\stdspace\theprimaryclass
\ifx\thesecondaryclass\relax\else; \thesecondaryclass\fi\par
{\bf Keywords}\stdspace \thekeywords\par}\vglue 7truept

}   

\ifplaintex
\font\phead=cmsl9 scaled 950
\font\pnum=cmbx10 scaled 913
\font\pfoot=cmsl9 scaled 950
\headline{\vbox to 0pt{\vskip -4.5mm\line{\small\phead\ifnum
\count0=\startpage ISSN 1472-2739 (on-line) 1472-2747 (printed)
\hfill {\pnum\folio}\else\ifodd\count0\def\\{ }%
\ifx\theshorttitle\relax\thetitle\else\theshorttitle\fi\hfill{\pnum\folio}
\else\def\\{ and }{\pnum\folio}\hfill\ifx\theshortauthors\relax\theauthors
\else\theshortauthors\fi\fi\fi}\vss}}
\footline{\vbox to 0pt{\vglue 0mm\line{\small\pfoot\ifnum\count0=\startpage
\copyright\ \gtp\hfill\else
\agt, Volume \thevolumenumber\ (\thevolumeyear)\hfill\fi}\vss}}
\else
\font=cmsl9 scaled 1050
\font\lnum=cmbx10 
\font=cmsl9 scaled 1050
\makeatletter
\def\@oddhead{{\small\lhead\ifnum\count0=\startpage ISSN 1472-2739 
(on-line) 1472-2747 (printed)\hfill {\lnum\number\count0}\else\ifodd\count0
\def\\{ }\ifx\theshorttitle\relax \thetitle \else\theshorttitle\fi\hfill
{\lnum\number\count0}\else\def\\{ and }{\lnum\number\count0}
\hfill\ifx\theshortauthors\relax 
\theauthors\else\theshortauthors\fi\fi\fi}}\def\@evenhead{\@oddhead}
\def\@oddfoot{\small\lfoot\ifnum\count0=\startpage\copyright\ \gtp\hfill\else
\agt, Volume \thevolumenumber\ (\thevolumeyear)\hfill\fi}
\def\@evenfoot{\@oddfoot}
\makeatother
\fi
\let\maketitlepage\makeagttitle
\let\makeshorttitle\maketitlepage
\let\maketitle\maketitlepage

\newwrite\gtoutfile
\long\gdef\makeheadfile{  
{\def\\{, }\def\s{ }
\immediate\openout\gtoutfile head.xxx
\immediate\write\gtoutfile{To: math@arxiv.org}
\immediate\write\gtoutfile{Subject: put OR rep NNNNN:ppppp}
\immediate\write\gtoutfile{--text follows this line--}
\immediate\write\gtoutfile{Proxy-for: \ifx\theasciiauthors\relax
\theauthors\else\theasciiauthors\fi\s<\ifx\theasciiemail\relax\theemail\else\theasciiemail\fi>}
\immediate\write\gtoutfile{\noexpand\\}
\immediate\write\gtoutfile{Authors: \ifx\theasciiauthors\relax
\theauthors\else\theasciiauthors\fi}
{\def\\{ }\immediate\write\gtoutfile{Title: \ifx\theasciititle\relax
\thetitle\else\theasciititle\fi}}
\immediate\write\gtoutfile{Subj-class: GT or SG, GR etc}
\immediate\write\gtoutfile{MSC-class: \theprimaryclass\ifx\thesecondaryclass\relax\else, \thesecondaryclass\fi}
\immediate\write\gtoutfile{Journal-ref: Algebr. Geom. Topol. \thevolumenumber\s
(\thevolumeyear) \startpage-\finishpage}
\immediate\write\gtoutfile{Comments: Published by Algebraic and
Geometric Topology at}
\immediate\write\gtoutfile{\s\s\s  http://www.maths.warwick.ac.uk/agt/AGTVol\thevolumenumber/agt-\thevolumenumber-\thepapernumber.abs.html}
\immediate\write\gtoutfile{\noexpand\\}
\immediate\write\gtoutfile{}
\ifx\theasciiabstract\relax
\immediate\write\gtoutfile{\theabstract}\else
\immediate\write\gtoutfile{\theasciiabstract}\fi
\immediate\write\gtoutfile{}
\immediate\write\gtoutfile{\noexpand\\}
\immediate\write\gtoutfile{}
\immediate\closeout\gtoutfile}}  

\def\maketitlepage{\makeagttitle\makeheadfile}
\let\makeshorttitle\maketitlepage
\let\maketitle\maketitlepage

\lognumber{3}
\volumenumber{3}
\volumeyear{2003}
\papernumber{3}
\published{28 January 2002}
\pagenumbers{89}{101}
\received{3 December 2002}
\accepted{23 January 2003}

\usepackage{amssymb, amsmath, epsf}
\def\hepsffile{\leavevmode\epsffile}
\nocolon
 
\theoremstyle{plain}
\newtheorem{thm}{Theorem}[subsection]

\newtheorem{prop}[thm]{Proposition}

\theoremstyle{definition}
\newtheorem{defin}[thm]{Definition}

\newtheorem{rem}[thm]{Remark}
\swapnumbers
\newtheorem{emf}[thm]{}
 
\def\pr{\protect\operatorname{pr}}
 
\def\Z{{\mathbb Z}}

\def\R{{\mathbb R}}
\def\N{{\mathbb N}}

\def\1{\hbox{\rm\rlap {1}\hskip.03in{\rom I}}}
\def\Bbbone{{\rm1\mathchoice{\kern-0.25em}{\kern-0.25em}
        {\kern-0.2em}{\kern-0.2em}I}}
 
\begin{document}
\hyphenation{Ca-m-po}
\title[The universal order one invariant of knots in $S^1$--fibered spaces]
{The universal order one invariant of framed knots\\in most $S^1$--bundles over orientable surfaces}
\asciititle{The universal order one invariant of framed knots in most S^1-bundles over orientable surfaces}

\author{Vladimir Chernov (Tchernov)}
\address{Department of Mathematics, 6188 Bradley Hall\\Dartmouth 
College, Hanover, NH 03755, USA}
\email{Vladimir.Chernov@dartmouth.edu}

\begin{abstract}
It is well-known that self-linking is the only $\Z$-valued 
Vassiliev invariant
of framed knots in $S^3$. However for most $3$-manifolds, in particular for
the total spaces of $S^1$-bundles over an orientable surface $F\neq S^2$, 
the space of $\Z$-valued order one invariants is infinite dimensional.  

We give an explicit formula for the order one invariant $I$ of framed 
knots in orientable 
total spaces of $S^1$-bundles over an orientable not necessarily compact 
surface
$F\neq S^2$. We show that if $F\neq S^2, S^1\times S^1,$ then
$I$ is the universal order one invariant, i.e.\ it
distinguishes every two framed knots that can be distinguished by order one
invariants with values in an Abelian group. 
\end{abstract}

\asciiabstract{It is well-known that self-linking is the only Z valued 
Vassiliev invariant of framed knots in S^3. However for most
3-manifolds, in particular for the total spaces of S^1-bundles over an
orientable surface F not S^2, the space of Z-valued order one
invariants is infinite dimensional.  We give an explicit formula for
the order one invariant I of framed knots in orientable total spaces
of S^1-bundles over an orientable not necessarily compact surface F
not S^2. We show that if F is not S^2 or S^1 X S^1, then I is the
universal order one invariant, i.e. it distinguishes every two framed
knots that can be distinguished by order one invariants with values in
an Abelian group.}


\maketitle

\section{Main Results}
\subsection{Introduction}
We work in the smooth category.

A {\em surface\/} is a not necessarily compact connected $2$-dimensional manifold.
A {\em curve\/} in a manifold is an immersion of $S^1$ into the manifold. A
{\em framed
curve\/} is a curve equipped with a transverse vector field. A {\em knot (framed
knot)\/} is an embedded curve. Knots and framed knots are studied up to ~the
corresponding isotopy equivalence relation. An $S^1$-bundle over a surface
is a locally-trivial $S^1$-fibration.

In~\cite{Chernovamsvolume} we used a new kind of skein relation 
to construct an order one invariant of unframed knots in an orientable 
total space $M$ of an $S^1$-bundle over a surface. 
The invariant takes value in a quotient of the
group ring of $H_1(M)$, and in the case of
spherical tangent bundles of surfaces it is 
a splitting of the Polyak's~\cite{Polyakbennequin} Arnold-Bennequin
type invariant of wave fronts on surfaces. 

However, as it is easy to see, the
invariant introduced in~\cite{Chernovamsvolume} is not universal i.e.\ 
there exist examples of two unframed knots that can 
not be distinguished by this
invariant, but can be distinguished by some other order one Vassiliev
invariant. 

In this paper we give an explicit geometric construction of 
the order one invariant $I$ of framed knots in an oriented total space $M$ 
of an
$S^1$-bundle over an oriented surface $F$, and we show that for $F\neq S^2,
S^1\times S^1$
the invariant $I$
distinguishes every two knots that can be distinguished by 
order one invariants with values in an Abelian group. 

The invariant $I$ takes values in the group of formal finite integer
linear combinations of the free homotopy classes of mappings of the wedge
of two circles into the manifold $M$, factorised by the automorphism of the
wedge that interchanges the circles. 
The geometric ideas we use to
construct $I$ are similar to those that allowed us~\cite{Chernovfronts} to obtain a formula for the
universal order one Arnold's~\cite{Arnoldsplit} 
$J^+$-type invariant of wave fronts on
an orientable surface $F\neq S^2$. 

In general first order invariants of knots and links in the total spaces of sphere-bundles 
appear to be very important in the study of wave 
propagation. 
In particular, recently the author and Yu.~Rudyak~\cite{ChernovRudyak} applied the order one invariants 
of links in the sphere-bundles to the study of the causality relation for wave fronts.

\subsection{Construction of the invariant}
Let $M$ be an oriented total space of an
$S^1$-bundle $p:M\rightarrow F$ 
over an oriented surface $F\neq S^2$. Clearly if two framed knots $K_1,
K_2\subset M$ belong to different components of the space of framed curves
in $M$, then they are not isotopic. For this reason when studying framed 
knots we restrict ourselves to a connected component $\mathcal F$ of the space of
framed curves in $M$. It is easy to verify that each connected component of the space of unframed curves in $M$ corresponds (under forgetting of the framing) to precisely two connected components of the space of framed curves. These two components are distinguished by the values of a spin structure on loops in the principal $SO(3)$-bundle  that are naturally associated with the framed curves. In turn, connected components of the space of unframed curves in $M$ are in the natural bijective correspondence with the conjugacy classes of the elements of $\pi_1(M)$.

Let $\mathcal C$ be a connected component of the space
of unframed curves in $M$ obtained by forgetting the framing on curves from
$\mathcal F$, and let $\mathcal L$ be the connected component of the space of
free loops on $F$ 
that contains the projections of curves from $\mathcal C$.
(Connected components of the space of free loops on $F$ are naturally
identified with the  conjugacy classes of the elements of $\pi_1(F)$.)

\begin{emf}\label{h-principle}{\bf $h$-principle for curves on $F$}\qua 
Clearly $\mathcal L$ contains many components of the space of curves
(immersions of $S^1$) on $F$. Put $\pr:STF\rightarrow F$ to be the
$S^1$-bundle obtained by the fiberwise spherization of the tangent bundle of
$F$. 
The $h$-principle~\cite{Gromov} says that the
space of curves on $F$ is weak homotopy equivalent to the space of free
loops $\Omega STF$ in $STF$. The equivalence is given by mapping a curve $C$
on $F$ to a loop $\vec C$ in $STF$ obtained by mapping $t\in S^1$ to the
point in $STF$ that corresponds to the velocity vector of $C$ at $C(t)$.
In particular, the connected components of the space of
curves on $F$ that are contained in $\mathcal L$ are naturally identified with
the connected components of the space of free loops in $STF$ that consist of
loops projecting to loops from $\mathcal L$.
\end{emf}

\begin{figure}[ht!]
 \begin{center}
  \epsfxsize 11cm
  \hepsffile{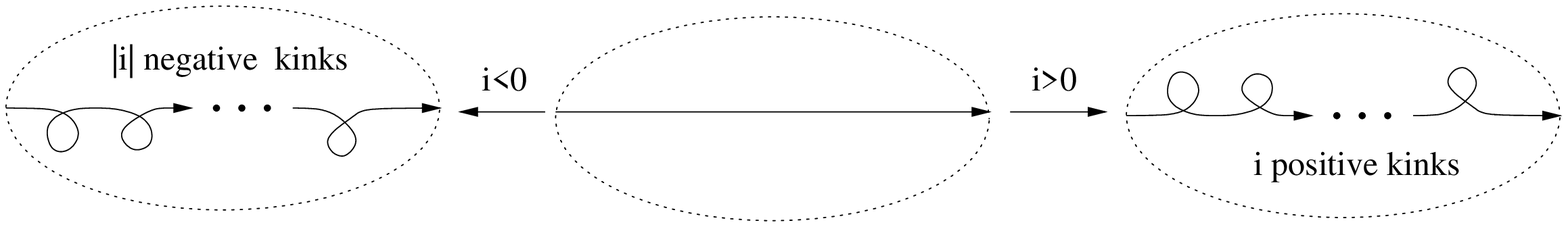}
 \end{center}
\caption{}
\label{kinks.fig}
\end{figure}

\begin{prop}\label{action}
Let $F\neq S^2$ be a (not necessarily compact) oriented surface, and let
$\mathcal C, \mathcal L$ be as above. 
Then the group $\Z$ acts freely and transitively on the set of connected
components of the space of curves on $F$ that realize loops from $\mathcal
L$. The action is as follows: $i\in\Z$ acts on a connected component $\mathcal K$ that
contains a curve $C\in\mathcal L$ by mapping it to the connected component 
$\mathcal K^i$ 
that contains the curve obtained from $C$ by the addition to it of $i$ positive
kinks, provided that $i\geq 0$; and by the addition of $|i|$ negative kinks
provided that $i<0$, see Figure~\ref{kinks.fig}.
\end{prop}
 
For the proof of Proposition~\ref{action} see Section~\ref{proofaction}.

We say that a curve $C\in\mathcal C$ is {\em generic\/} (with respect to the $S^1$-bundle
$p:M\rightarrow F$) if $p(C)$ is a curve and its only singularities are
double points of transverse self-intersection.

Fix a point $x_0\in M$.
Let $K_s\subset M$ be a singular knot whose only singularity is a transverse 
double point $d$. Let $B_2$ be the wedge of two oriented circles. Consider 
a mapping $\phi:S^1\rightarrow B_2$ that maps the preimages of the double 
point of $K_s$ to the base point of the wedge, respects the orientations
of the circles of the wedge, and is injective on the complement of the
preimage of the double point of $K_s$. Then there exists a unique mapping 
$\psi:B_2\rightarrow M$ such that $\psi \circ\phi =K_s$. Now we pull the
base point of $\psi(B_2)$ till it is located at the point $x_0$ and the two
loops of the wedge give an element of $\pi_1(M, x_0)\oplus\pi_1(M,x_0)$.

Clearly the free homotopy class of 
$\psi $ is well-defined modulo the action of
$\Z_2$ that interchanges the two circles. Hence to $K_s$ corresponds 
a unique element $b$
of the quotient set $\mathcal B$ of $\pi_1(M, x_0)\oplus\pi_1(M, x_0)$ modulo the actions of
$\pi_1(M)$ via conjugation (this action corresponds to the ambiguity in
choice of the path along which we pull the base point of the wedge till it is located at
$x_0$); 
and by the action of $\Z_2$ via permutation of the two
summands (this action corresponds to the ambiguity in the choice of the
first of the two loops of $K_s$).  
This $b\in \mathcal B$ is said to be {\em the element corresponding to the
singular knot $K_s$} (with one double point). 

Let $K$ be a knot in $M$ that is generic with respect to $p:M\rightarrow F$. 
Let $d$ be a double point of $p (K)$. Since $F$ is oriented we can
distinguish the two branches of $K$ over $d$.
We call the branch of $K$ over
$d$ {\em the left branch\/} if the $2$-frame that is formed by the
projections to $F$ of the velocity vectors of this branch and of the other branch
of $K$ over $d$ gives the chosen orientation of the oriented surface $F$.
The other branch is called the {\em right branch\/} of $K$ over $d$.

\begin{figure}[ht!]
 \begin{center}
  \epsfxsize 7cm
  \hepsffile{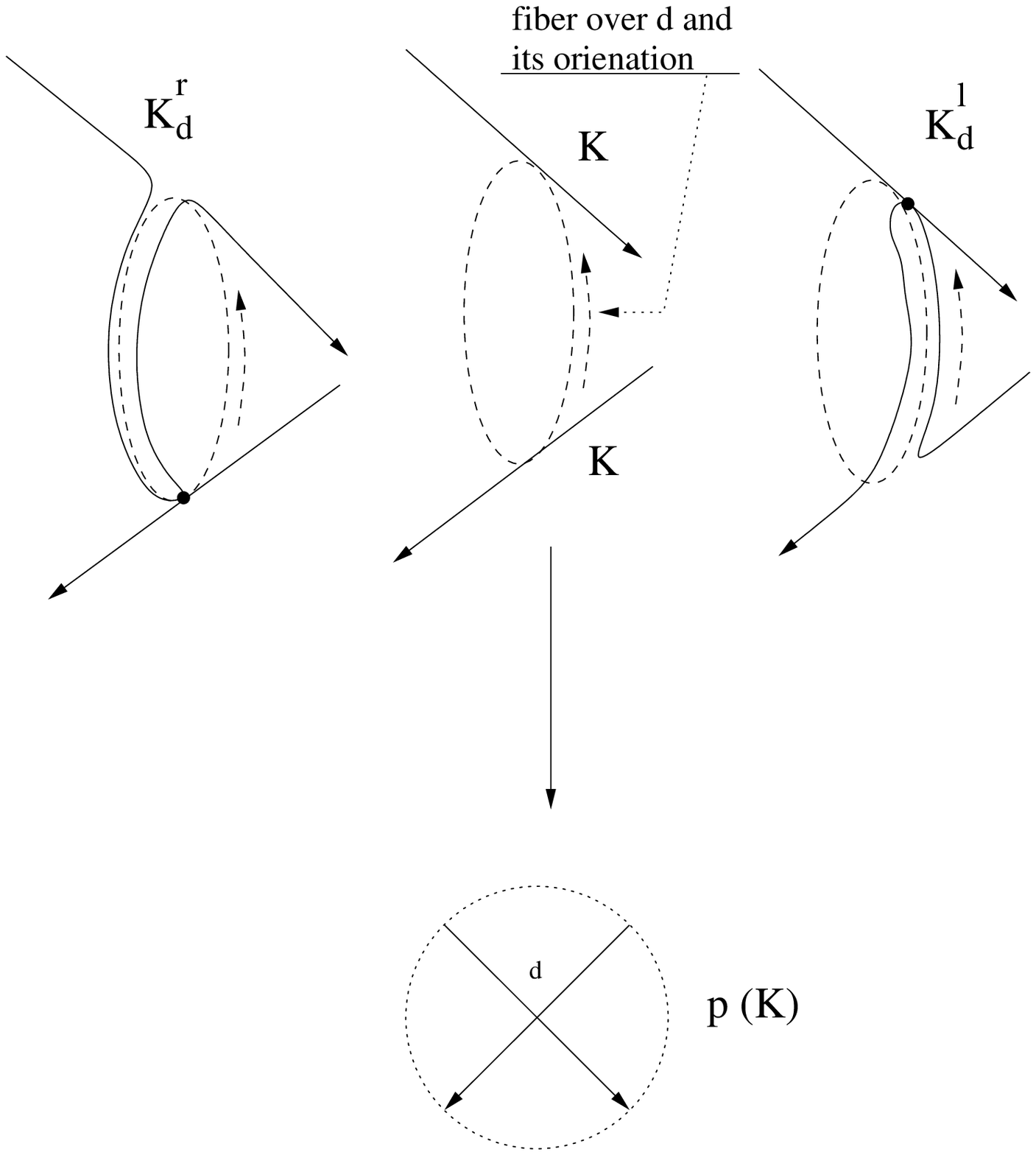}
 \end{center}
\label{model1.fig}
\caption{}
\end{figure}

Since both $M$ and $F$ are oriented, the $S^1$-fibers of $p:M\rightarrow F$
are naturally oriented.
To a double point $d$ of $p(K)$ we
associate two singular knots (with a double point) $K^r_d$ and $K^l_d$. 
The singular knots 
$K^r_d$ and $K^l_d$ are obtained from $K$ by taking respectively the right
and the left branch of $K$ over $d$ 
and pulling it along the
oriented $S^1$-fiber over $d$ in the direction coherent with the orientation
of the fiber till it intersects the other branch,
see Figure~\ref{model1.fig}.
Put $[K^r_d], [K^l_d]\in\mathcal B$ to be the elements that
correspond to $K^r_d$ and $K^l_d$.

Let $f\in\pi_1(M,x_0)$ be the class of the oriented $S^1$-fiber of
$p:M\rightarrow F$. As one can show, see~\ref{commute}, $f$ is in the center
of $\pi_1(M)$. Thus if $b\in\mathcal B$ is realizable as
$\alpha\oplus\beta\in\pi_1(M,x_0)\oplus\pi_1(M,x_0)$, then for $i,j\in\Z$ 
the element of
$\mathcal B$ realized by $\alpha f^i\oplus \beta f^j$ 
depends only on $b\in\mathcal B$ and $i,j\in\Z$ and does not depend on
the choice of realization of $b$ as $\alpha\oplus\beta$. We denote this
element of $\mathcal B$ by $[\alpha f^i, \beta f^j]$.
Put $[K,1]$ to be the element of $\mathcal B$ that corresponds to the free
homotopy class of the wedge with the base point at $K(1)$, the first loop
being $K$, and the second loop being trivial. 
Put $\Z[\mathcal B]$ to be the group
of formal finite linear combinations of the elements of $\mathcal B$ with integer
coefficients.

\begin{defin}[of the invariant $I$]\label{definI}
Fix a connected component $\mathcal K$ of the space of curves on $F$ that is
contained in the component $\mathcal
L$ of the space of free loops.

Let $K\in\mathcal F$ be a framed knot that is generic with respect to
$p:M\rightarrow F$ and has all the framing vectors non-tangent to the fibers
of $p$. For simplicity we assume that the projections of the framing vectors point to the right of $p(K)$. (Clearly every framed knot can be deformed to be the one with such properties 
by a $C^0$-small
deformation.) Put $i(K)$ to be the unique integer such that $p(K)\in\mathcal K^{i(K)}$,
see Proposition~\ref{action}. 

Define $I(K)\in\Z[\mathcal B]$ by
\begin{equation}
I(K)=i(K)([Kf,f^{-1}]-[Kf^{-1}, f])+
\sum_{d}2([K^l_d]-[K^r_d]).
\end{equation}
\end{defin}

\begin{defin}[Vassiliev invariants]
A transverse double point $t$ of a singular knot can be resolved in two 
essentially different ways. We say that a resolution of a double point is
{\em positive (resp. negative)\/} if the tangent vector to the
first strand, the tangent vector to the second strand, and the vector from
the second strand to the first form the positive $3$-frame (this does not depend
on the order of the strands).

A singular framed knot $K_s$ with $(n+1)$ 
transverse double points
admits $2^{(n+1)}$ possible resolutions of the double points. A sign of the resolution 
is put to be $+$ if the number of negatively resolved double points is even, and
it is put to be $-$ otherwise. 
Let $\mathcal A$ be an Abelian group and let $x$ be an $\mathcal A$-valued invariant of framed knots. The invariant $x$ is said to be of {\em finite
order\/} (or {\em Goussarov-Vassiliev invariant\/}) if there exists a positive 
integer $(n+1)$ such that for any singular knot $K_s$ with $(n+1)$
transverse double points the sum (with appropriate signs) of the values of $x$ on the nonsingular
knots obtained by the $2^{n+1}$ resolutions of the double points is zero. 
An invariant is said to be of order not greater than $n$ (of order $\leq n$) if $n$
can be chosen as integer in the above definition. The group of $\mathcal
A$-valued finite order invariants has an increasing filtration by the
subgroups of the invariants of order $\leq n$.
\end{defin}

\begin{thm}\label{order1}
Let $p:M\rightarrow F$ be an oriented $S^1$-bundle over a (not necessarily
compact) oriented surface $F\neq S^2$. Then
\begin{description}
\item[1] $I(K)$ is an isotopy invariant of the framed knot $K$;
\item[2] Let $K_s$ be a singular knot with one double point $d$, let $K^+_s$ and
$K^-_d$ be the nonsingular framed knots obtained by respectively the
positive and the negative resolution of $d$, and let
$\alpha,\beta\in\pi_1(M)$ be such that $[\alpha,\beta]$ is the element of
$\mathcal B$ that corresponds to $K_s$. Then
$I(K^+_s)-I(K^-_s)=2(-2[\alpha,\beta]+[\alpha f^{-1}, \beta f]+[\alpha f,
\beta f^{-1}])$, and
thus $I(K)$ is an order one invariant.
\end{description}
\end{thm}

For the Proof of Theorem~\ref{order1} see Section~\ref{prooforder1}.

The following Theorem says that $I(K)$ distinguishes all pairs of knots that can 
possibly be distinguished with the order one Vassiliev invariants, provided that 
$F\neq S^2, S^1\times S^1$. 
This means that $I$ is the universal order one invariant of knots in an
oriented total space $M$ of an $S^1$-bundle $p:M\rightarrow F$ over a
not necessarily compact oriented surface $F\neq S^2, S^1\times S^1$.

\begin{thm}\label{universal}
Let $M$ be an oriented $3$-manifold, let $F\neq S^2, S^1\times S^1$ be an
oriented (not necessarily compact) surface, and let $p:M\rightarrow F$ be an
$S^1$-bundle. 
Let $\mathcal F$ be a connected component of the
space of framed curves in $M$, and let $K_1, K_2\in\mathcal F$ be framed
knots. 
Let $\widetilde I$ be an order one Vassiliev invariant (with values in some Abelian group)  
such that $\widetilde I (K_1)\neq \widetilde I(K_2)$, then $I(K_1)\neq I(K_2)$.
\end{thm} 

For the proof of Theorem~\ref{universal} see Section~\ref{proofuniversal}.

\begin{rem}The statement of Theorem~\ref{universal} holds also
in the case of $p:S^1\times S^1\times S^1=ST(S^1\times S^1)\rightarrow S^1\times S^1$. The
proof of the Theorem for this case is obtained by a straightforward
generalization.
\end{rem}

\section{Proofs}
\subsection{Proof of Proposition~\ref{action}}\label{proofaction}
We start with the following Propositions.

\begin{prop}\label{commute}
Let $N, L$ be oriented
manifolds and let $q:N\rightarrow L$ be an $S^1$-bundle.  
Then the 
class $f\in\pi_1(N)$ of the oriented $S^1$-fiber of $q$ 
is in the center of $\pi_1(N)$.
\end{prop}

Take $\alpha\in\pi_1(N)$. Consider $\mu:S^1\times S^1\rightarrow N$ with
$\mu\big|_{t\times S^1}$ being the oriented $S^1$-fiber of $q$ that contains
$\alpha(t)$. Then the restriction of $\mu$ to the $2$-cell of the torus 
gives the commutation relation between $\alpha$ and $f\in\pi_1(N)$.\qed

\begin{prop}[A.~Preissman]\label{Preissman}
Let $F\neq S^2, S^1\times S^1$
be an oriented  (not necessarily compact) surface 
and let $G$ be a nontrivial commutative subgroup of $\pi_1(F)$.
Then $G$ is infinite cyclic.
\end{prop}

\begin{emf}{\bf Proof of Proposition~\ref{Preissman}}

It is well known that any closed oriented $F$, other than $S^2, S^1\times
S^1$,
admits a hyperbolic metric of a constant negative curvature.
(It is induced from the universal covering of $F$ by the hyperbolic plane
$H$.)
The Theorem by A.~Preissman (see~\cite{Docarmo} pp.\ 258--265)
says that if $M$ is a compact Riemannian manifold with a negative curvature,
then any nontrivial Abelian subgroup $G<\pi_1(M)$ is infinite cyclic.

If $F$ is not closed, then the statement of the Proposition is
true, since $\pi_1(F)$ is a free group on a countable or finite set of generators, see 
Ahlfors and Sario~\cite{AhlforsandSario}, chapter 1, or~\cite{Massey}, pp.\ 143 and 199--200.
\qed
\end{emf}

\begin{emf}
The proof is based on the $h$-principle, see~\ref{h-principle}.
Let $f\in\pi_1(STF)$ be the class of the oriented $S^1$-fiber of $\pr:STF\rightarrow
F$. Proposition~\ref{commute} says that $f$ is in the center of $\pi_1(STF)$.

Let $C$ be a curve from $\mathcal L$. Take $i\in\Z$ and put $C^i$ to be a
curve obtained from $C$ by the addition of $i$ positive kinks for $i\geq 0$; 
and by the addition of $|i|$ negative kinks for $i<0$. It is easy to see that 
$\vec {C^i}=\vec C
f^i\in\pi_1(STF)$. Since $\ker (\pr_*:\pi_1(STF)\rightarrow \pi_1(F))$ is
generated by $f$, the $h$-principle implies that the action of $\Z$
(introduced in Proposition~\ref{action}) on the
set of connected components of the space of curves that are contained in
$\mathcal L$ is well defined and transitive.

To show that the action is free it suffices to show that for any $\alpha,
\beta\in\pi_1(STF)$ if 
\begin{equation}\label{conjugation}
\alpha\beta\alpha^{-1}=\beta f^k,
\end{equation}
then $k=0$.

If $F=T^2=S^1\times S^1$, then $STT^2=S^1\times S^1\times S^1$ and the fact
that $k=0$ is obvious. For this reason below we assume that $F\neq T^2$. 
Clearly for such $\alpha, \beta$ the elements $\pr_*(\alpha)$ and $\pr_*(\beta)$
commute in $\pi_1(F)$. Proposition~\ref{Preissman} implies that there exist
$\bar g\in\pi_1(F)$ and $i,j\in\Z$ such that $\pr_*(\alpha)=\bar g^i$ and 
$\pr_*(\beta)=\bar g^j$. Take $g\in\pi_1(STF)$ such that $\pr_*(g)=\bar g$. 
Since $f$ is in the center of $\pi_1(STF)$ and $f$ generates $\ker \pr_*$, we
get that there exist $l,m\in\Z$ such that $\alpha=g^if^l$ and
$\beta=g^jf^m$. Substitute these expressions for $\alpha$ and $\beta$
into~\eqref{conjugation} and use the fact that $f$ is in the center of
$\pi_1(STF)$ to get that $f^k=1$. Since $\pi_2(STF)=0$ for our manifolds
$STF$, we
get that $f$ has infinite order in $\pi_1(STF)$. Thus $k=0$. This finishes the
proof of Proposition~\ref{action}.\qed
\end{emf}

\subsection{Proof of Theorem~\ref{order1}}\label{prooforder1}
Let $K_1, K_2$ be two isotopic oriented framed knots such that $p(K_1), p(K_2)$ are immersions, the framing vectors of knots 
are nowhere tangent to
the fibers of $p:M\rightarrow F$ and project to the nonzero vectors pointing to the right of $p(K_1)$ and $p(K_2)$, respectively.

\begin{figure}[ht!]
 \begin{center}
  \epsfxsize 11cm
  \hepsffile{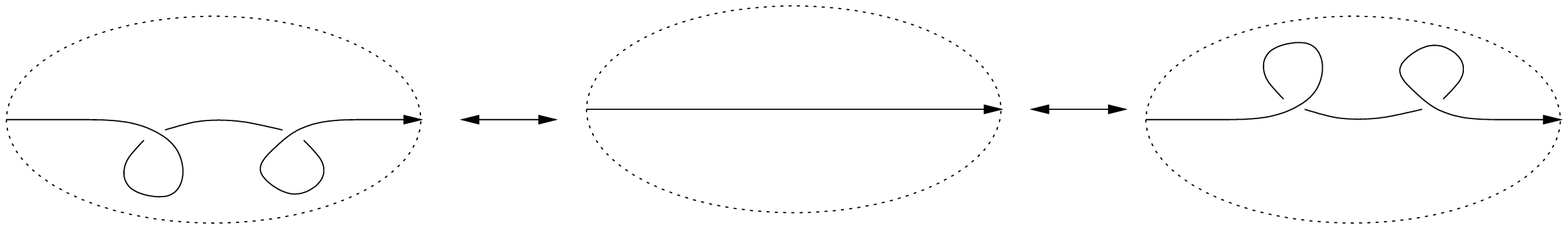}
 \end{center}
\caption{}
\label{move1.fig}
\end{figure}

Then it is clear that there is an isotopy between $K_1$ and $K_2$ 
that can be decomposed 
into 
\begin{description}\item[1]
isotopies that project to the ambient isotopies of
projections with the framing vectors nowhere tangent to the fibers of $p$ and projections of them pointing to 
the right from the oriented knot projections; and 
\item[2] the sequence of moves such that 
\begin{description} 
\item[a] they happen in the charts of $M$
homeomorphic to $\R^3=(x,y,z)$ with the lines $(x_0, y_0, z)$ 
for fixed $(x_0, y_0)$
being the arcs of the $S^1$-fibers of $p$; 
\item[b] projections of the moves to the $(x,y)$-plane correspond to the
second and third Reidemeister moves, and the first Reidemeister move for
framed knots with blackboard framing shown in Figure~\ref{move1.fig} and
its reflections. (At the start and end of these moves the framing vectors are assumed to be nowhere tangent to the fibers of 
$p$ and their projections point to the right from the oriented knot projections.)
\end{description}
\end{description}

The invariance of $I$ under the isotopies that project to the ambient
isotopies of projections is obvious.

The increments into $I$ that correspond to the double points of $p(K)$ that do 
not participate in the Reidemeister moves are unchanged under the moves.

Clearly the connected component of the space of curves on $F$ that contains
$p(K)$ is unchanged under the second and third Reidemeister
moves.

A straightforward verification shows that the summands
$2([K^l_{d_1}]-[K^r_{d_1}])$ and $2([K^l_{d_2}]-[K^r_{d_2}])$ corresponding
to the two extra double points $d_1$ and $d_2$ of $p(K)$ that participate in an oriented version of
the second Reidemeister move cancel out, and thus $I$ is invariant under the
second move.

There is a natural correspondence between the three branches of $K$ that are
present on the diagram before and after the third move. This correspondence
induces the natural identification between the three double points of $p(K)$
before the move and after the move. (We identify the two points that are
the double points of the projection of the corresponding pairs of branches.) 
Now it is easy to see that for the corresponding double points
$d$ and $d'$ of $p(K)$ the summands $2([K^l_{d}]-[K^r_{d}])$ and
$2([K^l_{d'}]-[K^r_{d'}])$ are equal.
Thus $I$ is invariant under the third move.

Clearly the number $i(K)$ is changed by $\pm 2$ under the first move 
(depending on the version of the move that takes place). 
Thus the summand $i(K)([Kf,f^{-1}]-[Kf^{-1}, f])$ increases by $\pm
2([Kf,f^{-1}]-[Kf^{-1}, f])$. On the other hand it is easy to see that the
increments into $\sum_{d}2([K^l_d]-[K^r_d])$ corresponding to the two
double points of $p(K)$ that participate in the first move also do not
cancel out and the sum $\sum_{d}2([K^l_d]-[K^r_d])$ increases by $\mp
2([Kf,f^{-1}]-[Kf^{-1}, f)$. Thus $I(K)=i(K)([Kf,f^{-1}]-[Kf^{-1}, f])+
\sum_{d}2([K^l_d]-[K^r_d])$ does not change under the first move.
This shows that $I(K)$ is invariant under isotopy.

The proof of the second statement of the Theorem is a straightforward
calculation.\qed

\subsection{Proof of Theorem~\ref{universal}}\label{proofuniversal}
Let $\mathcal A$ be an Abelian group and let 
$\widetilde I$ be an $\mathcal A$-valued 
order one Vassiliev invariant of framed knots from $\mathcal F$.
Since $\widetilde I$ is an order one invariant, it is defined (up to an additive constant) by its derivative $\widetilde I'$ 
i.e.\ by the values of its
increments under the passages through the codimension zero strata of the discriminant subspace of $\mathcal
F$. (The discriminant is the subspace of $\mathcal F$ formed by
singular knots, and the codimension zero strata of the discriminant are formed by singular knots whose only singularity is one transverse double point.) Since $\widetilde I$ is an order one invariant, we get that the values of
the increments depend only on the elements of $\mathcal B$ that correspond to
singular knots with one double point that we obtain when we cross the
discriminant.

Since we have fixed the connected component $\mathcal C$, we get that for any $[\alpha_1, \alpha_2]\in\mathcal B$ 
that corresponds to a singular knot from $\mathcal F$ the loop
$\alpha_1\alpha_2$ is free homotopic to a curve from $\mathcal C$.
Observe that for all $[\alpha_1, \alpha_2]$ that participate in the
definition of $I(K)$ the loop $\alpha_1\alpha_2$ is also free homotopic to 
a loop from $\mathcal C$.

{\em For this reason by abuse of notation below in the proof 
we denote by $\mathcal B$ 
the subset of $\mathcal B$ that consists of the elements realizable by  
$\alpha_1\oplus\alpha_2\in\pi_1(M)\oplus\pi_1(M)$ with the loop
$\alpha_1\alpha_2$ free homotopic to curves from $\mathcal C$.\/}

Consider the homomorphism
$g:\Z[\mathcal B]\rightarrow\Z[\mathcal B]$
that maps 
\begin{equation}
[s_1, s_2]\rightarrow2(-2[s_1, s_2]+[s_1f,
s_2f^{-1}]+[s_1f^{-1}, s_2f]).
\end{equation}
(This homomorphism describes
the behavior of $I$ under crossings of the discriminant, see
Theorem~\ref{order1}.2. Recall that
by~\ref{commute} $f$ is in the center of $\pi_1(M)$.)
To prove the Theorem it suffices to show that $\ker g=0$.

Let $\overline {\mathcal B}$ be the quotient set of $\pi_1(F)\oplus\pi_1(F)$
via the actions of $\pi_1(F)$ that acts by conjugation of both summands and
by the action of $\Z_2$ that acts by permuting the summands. Once again by abuse of notation below we denote by $\overline {\mathcal B}$ the part of $\overline {\mathcal B}$ that is formed by the classes of $\alpha_1\oplus \alpha_2 \in \pi_1(F)\oplus\pi_1(F)$ such that $\alpha_1\alpha_2$ is free homotopic to the loops from $\mathcal L$.
Let $q:\mathcal B\rightarrow \overline{\mathcal B}$ be the natural mapping induced by $p_*:\pi_1(M)\rightarrow \pi_1(F)$. (One verifies that this mapping is really well-defined.)

One verifies that $\Z[\mathcal B]$ splits into the direct sum over
$\overline{\mathcal B}$ of $\Z$-submodules
that are finite linear combinations of the elements of $\mathcal B$
projecting to the same element of $\overline {\mathcal B}$. Clearly $g$ maps every
summand to itself. Thus it suffices
to show that the restriction of $g$ to every summand has trivial kernel.

Fix $\bar b\in \overline{\mathcal B}$. Below we construct the ordering on
$q^{-1}(\bar b)$, that makes it isomorphic (as an ordered set)
to $\N$ or to $\Z$ (depending on $\bar b$). One verifies that the matrix
of the restriction of $g$ to $\Z[q^{-1}(\bar b)]$
written with respect to the basis that is the
ordered set $q^{-1}(\bar b)$ appears to be tridiagonal with all nonzero entries on the
diagonal below the main one. Thus the restriction of $g$ to 
$\Z[q^{-1}(\bar b)]$ has trivial kernel, and this 
proves the Theorem.

To construct the ordering on $q^{-1}(\bar b)$ we need the following
proposition.

\begin{prop}\label{technical}
Let $F\neq S^2, S^1\times S^1$ be a (not necessarily compact) oriented surface, 
let $p:M\rightarrow F$ be an
$S^1$-bundle with oriented $M$, let $f\in\pi_1(M)$ be the class of the
oriented fiber of $p$, and let $\alpha_1, \alpha_2$ be elements
of $\pi_1(M)$.
\begin{description}
\item[a] $\alpha_1$ and $\alpha_2$ commute in $\pi_1(M)$ if and only if
$p_*(\alpha_1)$ and $p_*(\alpha_2)$ commute in $\pi_1(F)$.
\item[b] If $p_*(\alpha_1)$ and $p_*(\alpha_2)$ are conjugate in
$\pi_1(F)$, then there exists a unique $i\in\Z$ such that $\alpha_1$ and
$\alpha_2 f^i$ are conjugate in $\pi_1(M)$.
\item[c]
Let $\beta_1, \beta_2\in\pi_1(M)$ be such that
$(\delta \alpha_1 \delta^{-1}, \delta \alpha_2 \delta^{-1})=(\beta_1,
\beta_2)\in\pi_1(M)\oplus\pi_1(M)$, for some $\delta\in\pi_1(M)$.
If there exists $\xi\in\pi_1(F)$ such that 
$p_*(\alpha_2)=\xi p_*(\alpha_1) \xi ^{-1}$ and $p_*(\alpha_1)=\xi p_*(\alpha_2)
\xi^{-1}$, then
$p_*(\alpha_1)=p_*(\alpha_2)$, $p_*(\beta_1)=p_*(\beta_2)$;
and hence there exist unique $i,j\in\Z$ such that
$\alpha_1=\alpha_2 f^i$, $\beta_1=\beta_2f^j$. Moreover $i=j$.
\end{description}
\end{prop}

The proof of the proposition is a straightforward calculation (similar to the
one we did when proving that the action of $\Z$ introduced in~\ref{action} is
free) and is based on
Propositions~\ref{commute}, \ref{Preissman}, and the fact that
$f$ generates $\ker p_*$ and has infinite order.

{\em The ordering of the basis $q^{-1}(\bar b)$ of $\Z[q^{-1}(\bar b)]$ such that the matrix
of $g\big|_{\Z[q^{-1}(\bar b)]}$ written with respect to this ordered basis
is tridiagonal with all the elements on the diagonal below the main one being
nonzero is constructed as follows:\/}

a)\qua If $b\in q^{-1}(\bar b)$ can be realized as $(\alpha_1, \alpha_2)$ such that
$\xi p_*(\alpha_1) \xi^{-1}=p_*(\alpha_2)$ and
$\xi p_*(\alpha_2) \xi^{-1}=p_*(\alpha_1)$,
for some $\xi\in\pi_1(F)$,
then any realization of any element of $q^{-1}(\bar b)$ has this property.
From~\ref{technical}.c we get that
every element $b\in q^{-1}(\bar b)$ determines
a unique $i\in \N$ such that $b$ can be realized as $(\alpha_1,
\alpha_2)$ with $\alpha_1 f^i=\alpha_2$. One verifies that these
natural numbers are different for different elements of $q^{-1}(\bar b)$.
(Recall that as it was said in the beginning of the proof of the Theorem,
$\mathcal B$ in this proof denotes the subset of the original $\mathcal B$
that consists of elements realizable by $\alpha_1\oplus\alpha_2$ with
$\alpha_1\alpha_2$ being a loop free homotopic to curves from the fixed 
connected component $\mathcal C$ of the space of curves in $M$.)
The ordering on $q^{-1}(\bar b)$ is induced by the magnitude of $i\in\N$
and it makes $q^{-1}(\bar b)$ isomorphic to $\N$ as the ordered set.

b)\qua If $b\in q^{-1}(\bar b)$ can not be realized as an element of the type described above,
then none of the elements of  $q^{-1}(\bar b)$ can. This allows us to
distinguish one loop of $\bar b$, and consequently to distinguish one loop 
of the elements of $q^{-1}(\bar b)$. We use the $\Z_2$ action on
$\pi_1(M)\oplus\pi_1(M)$ (used to introduce $\mathcal B$) to interchange the
two loops, so that the first loop projects to the distinguished loop of
$\bar b$. We get that
every element of $q^{-1}(\bar b)$
can be realized in a unique way as an element of
the set $\widetilde{\mathcal B}$ that is the quotient of
$\pi_1(M)\oplus\pi_1(M)$ modulo the action of $\pi_1(M)$ by
conjugation of both summands.
If $(s_1, s_2)$ and $(s_3, s_4)\in \widetilde{\mathcal B}$ realize two elements of
$q^{-1}(\bar b)$, then
there exists a unique $i\in\Z$ such that $s_1 f^i$ is conjugate to $s_3$,
see~\ref{technical}.b. As it was said in the beginning of the proof,
$s_1s_2$ and $s_3s_4$ are conjugate in $\pi_1(M)$, since they correspond
to knots from the same connected component $\mathcal C$ of the space of
curves in $M$.
One uses this to verify that
if $i=0$, then $(s_1, s_2)$ and $(s_3, s_4)$ realize the same
element of $q^{-1}(\bar b)$. The
ordering on $q^{-1}(\bar b)$ is induced by the magnitude of $i$, and it makes
$q^{-1}(\bar b)$ isomorphic to $\Z$ as the ordered set.

This finishes the proof of Theorem~\ref{universal}. \qed

\rk{Acknowledgments}

This work was supported by the free term research money from Dartmouth College and the author is thankful to Dartmouth College for the support.

\Addresses

\end{document}